\documentclass[preprint,12pt]{elsarticle}

\usepackage{amssymb}
\usepackage{mathrsfs}
\usepackage{stmaryrd}
\input diagxy

\newtheorem{Th}{\bf Theorem}[section]
\newtheorem{Lem}[Th]{\bf Lemma}
\newtheorem{Pro}[Th]{\bf Proposition}

\newtheorem{Def}[Th]{\bf Definition}
\newtheorem{Exam}[Th]{\bf Example}

\usepackage{amssymb}

\newcommand{\p}[1]{{\bf Proof.} #1 \ $\Box$}

\journal{Topology and Its Applications}

\begin{document}

\begin{frontmatter}
\title{Algebraic representation of continuous lattices via the open filter monad, revisited}
\author{Wei YAO}
\address{School of Mathematics and Statistics, Nanjing University of Information Science and Technology, Nanjing, China}
\author{Yueli Yue}

\address{School of Mathematical Sciences, Ocean University of China, Qingdao, China}
\end{frontmatter}

In [A. Day, Filter monads, continuous lattices and closure systems,
Can. J. Math. 27 (1975) 50--59], Day showed that continuous lattices are precisely the algebras of the open filter monad over the category of $T_0$ spaces.

The aim of this paper is to give a clean  and clear version of the whole process of Day's approach.

\section{Preliminaries}

\subsection{Monad and the related algebra}

Let $F,G : {\bf A}\longrightarrow{\bf  B}$ be two functors. A {\it natural transformation} $\tau$
from $F$ to $G$ (denoted by $\tau: F\longrightarrow  G$) is a function that assigns
to each {\bf A}-object $A$ a {\bf B}-morphism $\tau_A : F(A)\longrightarrow G(A)$ such that $G(f)\cdot\tau_A=\tau_{A'}\cdot F(f)$ holds for each
{\bf A}-morphism $f:A\longrightarrow A'$.

A {\it monad} over a category {\bf X} is a triple $(T,\eta,\mu)$ consisting of
a functor $T:{\bf X}\longrightarrow{\bf X}$ and two natural transformations
$\eta:id_{\bf X}\longrightarrow  T$ and $\mu:T\circ T\longrightarrow T$
such that $\mu\cdot T\mu=\mu\cdot \mu T$ and $\mu\cdot\eta T=id_T=\mu\cdot T\eta$.

Given a monad $(T,\eta,\mu)$ over {\bf X}, a $T$-algebra (or an Eilenberg-Moore algebra) is a pair $(X,r)$,
where $X$ is an {\bf X}-object and $r:T(X)\longrightarrow X$ is a structured morphism satisfying that $r\cdot T(r)=r\cdot \mu_X$ and $id_X=r\cdot \eta_X$.

\subsection{The open filter monad over the category of $T_0$ topological spaces}

\begin{Def}Let $(X,\mathcal{O}(X))$ be a topological space.
A subset $v\subseteq\mathcal{O}(X)$ is called
an open filter of $X$ if $v$ is a lattice-theoretic filter of $(\mathcal{O}(X),\subseteq)$, that is,

{\rm(F)} $A\cap B\in u$ if{}f $A,B\subseteq u\ (\forall A,B\in\mathcal{O}(X))$.
\end{Def}

Let $\Phi(X)$ denote the set of all open filters of $X$. Then the pair $(\Phi(X),\subseteq)$ is a dcpo, where for each directed subset $\mathcal{A}$ of $\Phi(X)$, $\bigvee\mathcal{A}=\bigcup\mathcal{A}$.

\begin{Exam}Let $(X,\mathcal{O}(X))$ be a topological space. For $A\in 2^X$, define $[A]\subseteq\mathcal{O}(X)$ by $[A]=\{B\in\mathcal{O}(X)|\ A\subseteq B\}$. Then $[A]\in\Phi (X)$. For every $x\in X$, denote $[\{x\}]$ by $[x]$. It is easy to show that,

{\rm(1)} for every $v\in\Phi (X)$, $v=\bigcup\limits_{A\in v}[A]$;

{\rm(2)} for each $B\in\mathcal{O}(X)$, $B\in v$ if{}f $[B]\subseteq v$.
\end{Exam}

Let $X$ be a topological space. For every $A\in\mathcal{O}(X)$, define $\phi(A)\subseteq \Phi(X)$ by $\phi(A)=\{v\in\Phi(X)|\ A\in v\}$. Then $\{\phi(A)|\ A\in\mathcal{O}(X)\}$ becomes a topological base, which generates a $T_0$ topology, denoted by $\mathcal{O}(\Phi(X))$. It is easy to see that every member of $\mathcal{O}(\Phi(X))$ is an upper set in the dcpo $(\Phi(X),\subseteq)$.

Let $f:X\longrightarrow Y$ be a mapping. Define $\Phi(f):\Phi(X)\longrightarrow \Phi(Y)$ by $$\Phi(f)(u)=\{B\in\mathcal{O}(Y)|\ f^\leftarrow(B)\in u\}.$$It is easily shown that $\Phi(f)$ is a well-defined mapping.
Then $\Phi$ is an endofunctor on ${\bf Top}_0$.

Define $\eta_X:X\longrightarrow\Phi(X)$ and $\mu_X:\Phi^2(X)\longrightarrow\Phi(X)$ respectively by $\eta_X(x)=[x]\ (\forall x\in X)$ and $\mu_X(\alpha)(A)=\alpha(\phi(A))\ (\forall \alpha\in\Phi^2(X),\forall A\in\mathcal{O}(X)))$.

\begin{Pro}Both $\eta:id_{{\bf Top}_0}\longrightarrow \Phi $ and $\mu:\Phi^2\longrightarrow\Phi$ are natural transformations, and further $(\Phi,\eta,\mu)$ is a monad over ${\bf Top}_0$.\end{Pro}

\subsection{Continuous lattices and the Scott topology}

\begin{Def}{\rm (1)} A nonempty subset $I$ of a poset $X$ is
called an ideal if it is a lower directed subset. Denote by $\mathcal{I}(X)$ the set of all ideals of $X$.
A poset $X$ is called a dcpo if $\bigvee I$ exists for every $I\in\mathcal{I}(X)$.

{\rm (2)} Let $X$ be a dcpo and $x,y\in X$.
If for every $I\in\mathcal{I}(X)$, $y\leq\bigvee I$ always implies $x\in I$,
then $x$ is called way-below $y$, in symbols, $x\ll y$. For each $x\in X$, denote ${\Uparrow}x=\{y\in X|\ x\ll y\}$ and ${\Downarrow}x=\{y\in X|\ y\ll x\}$.

{\rm(3)} A complete lattice $X$ is called a continuous lattice if $x=\bigvee{\Downarrow}x$ holds for every $x\in X$.\end{Def}

\begin{Def}On a dcpo $X$, a subset $V\subseteq X$ is called {\it Scott open} if $V$ is an upper set and $\bigvee I\in V$ implies $I\cap V\not=\emptyset$ for every $I\in\mathcal{I}(X)$. The set of all Scott open subsets of $X$ forms a topology, called the {\it Scott topology} on $X$, denoted by $\sigma(X)$.\end{Def}

\begin{Pro}In a continuous lattice $X$, it holds that

{\rm(1)} $\forall x,y\in X$, if $x\ll y$, then there exists $z\in X$ such that $x\ll z\ll y$;

{\rm(2)} the family $\{{\Uparrow}x|\ x\in X\}$ forms a base of $\sigma(X)$.\end{Pro}

\begin{Pro}\label{wbI} Let $X$ be a dcpo, $I\in\mathcal{I}(X)$ and $x\in X$. Then

{\rm(1)} If $x\leq \bigvee I$, then ${\Downarrow}x\subseteq I$;

{\rm(2)} If $x\ll \bigvee I$, then $x\in I$.\end{Pro}

Let $X$ be a complete lattice equipped with the Scott topology $\sigma(X)$. For every $v\in\Phi(X)$,
define
$v^{\downarrow}=\bigcup_{A\in v}A^{\downarrow},$ where $A^{\downarrow}=\{y\in X|\ A\subseteq{\uparrow}y\}$.
It is a routine to show that $v^{\downarrow}\in\mathcal{I}(X)$ for every $v\in\Phi(X)$.

\begin{Pro}
Let $X$ be a complete lattice. Then $X$ is a continuous lattice if{}f $x=\bigvee [x]^{\downarrow}$ holds for every $x\in X$.\end{Pro}

\p{$\Longrightarrow$: We only need to show that ${\Downarrow}x\subseteq[x]^{\downarrow}\subseteq{\downarrow}x$. Firstly, the inclusion $[x]^{\downarrow}\subseteq{\downarrow}x$ is obvious. Secondly, for every $y\ll x$, we have $x\in{\Uparrow}y\in\sigma(X)$ and then ${\Uparrow}y\in[x]$. Thus $y\in({\uparrow}y)^{\downarrow}\subseteq({\Uparrow}y)^{\downarrow}\subseteq[x]^{\downarrow}$. By the arbitrariness $y$, ${\Downarrow}x\subseteq[x]^{\downarrow}$, as desired.

$\Longleftarrow$: We only need to show that $[x]^{\downarrow}\subseteq{\Downarrow}x$. If $y\in[x]^{\downarrow}$, then there exists $A\in\sigma(X)$ such that $y\in A^{\downarrow}$, and then $y\leq x$. Let $I\in\mathcal{I}(X)$ and $x\leq\bigvee I$. Then $\bigvee I\in A$ and so there exist $z\in A\cap I$. Thus $y\leq z$ and then $y\in I$. Hence, $y\ll x$, and consequently, $[x]^{\downarrow}\subseteq{\Downarrow}x$.}

\section{The representation of continuous lattices as monad algebras}

\begin{Th}Let $X$ be a continuous lattice equipped with the Scott topology and define $r:\Phi(X)\longrightarrow X$ by $r(v)=\bigvee v^{\downarrow}$. Then the pair $(X,r)$ is a $\Phi$-algebra over the monad $(\Phi,\mu,\eta)$.\end{Th}

\begin{Lem}For every $A\in \sigma(X)$, it holds that $r^\leftarrow(A)\subseteq \phi(A)$.\end{Lem}

\p{If $v\in r^\leftarrow(A)$, then $\bigvee v^{\downarrow}=r(v)\in A$ and $v^{\downarrow}\cap A\not=\emptyset$ since $v^{\downarrow}$ is an ideal. Thus there exists $x\in A$ such that $x\in v^{\downarrow}$. For $x\in v^{\downarrow}$, there exists $B\in v$ such that $x\in B^{\downarrow}$. Then $B\subseteq{\uparrow}x\subseteq A$. Therefore, $A\in v$ and $v\in\phi(A)$. Hence, $r^\leftarrow(A)\subseteq \phi(A)$.}

\begin{Lem}For every $A\in\sigma(X)$ and every $x\in X$, if $x\ll \bigwedge A$, then $\phi(A)\subseteq r^\leftarrow({\Uparrow}x).$\end{Lem}

\p{If $v\in\phi(A)$, then $A\in v$. Then $\bigwedge A\in v^{\downarrow}$ and $\bigwedge A\leq\bigvee v^{\downarrow}=r(v)$.
If $x\ll \bigwedge A$, then $x\ll r(v)$ and $v\in r^\leftarrow({\Uparrow}x)$. By the arbitrariness of $v$, we have $\phi(A)\subseteq r^\leftarrow({\Uparrow}x)$.}

\begin{Pro}$r\cdot\Phi(r)=r\cdot \mu_X$.\end{Pro}

\p{Let $\alpha\in\Phi^2(X)$. Then
$$r(\Phi(r)(\alpha))=\bigvee(\Phi(r)(\alpha))^{\downarrow},\ \ \ r(\mu_X(\alpha))=\bigvee(\mu_X(\alpha))^{\downarrow}.$$

Firstly, $$\begin{array}{lll}(\Phi(r)(\alpha))^{\downarrow}&=&\{x\in X|\ \exists A\in\sigma(X)\ s.t.\ A\in\Phi(r)(\alpha)\ and\ A\subseteq{\uparrow}x\}
\\&=&\{x\in X|\ \exists A\in\sigma(X)\ s.t.\ r^\leftarrow(A)\in\alpha\ and\ A\subseteq{\uparrow}x\},\end{array}$$
and
$$\begin{array}{lll}(\mu_X(\alpha))^{\downarrow}&=&\{x\in X|\ \exists A\in\sigma(X)\ s.t.\ A\in\mu_X(\alpha)\ and\ A\subseteq{\uparrow}x\}
\\&=&\{x\in X|\ \exists A\in\sigma(X)\ s.t.\ \phi(A)\in\alpha\ and\ A\subseteq{\uparrow}x\}.\end{array}$$

Firstly, by Lemma 2.2, $(\Phi(r)(\alpha))^{\downarrow}\subseteq(\mu_X(\alpha))^{\downarrow}$ and then $r(\Phi(r)(\alpha))\leq r(\mu_X(\alpha)).$

Secondly, let $x_0=r[\mu_X(\alpha)]$. Then $x_0=\bigvee{\Downarrow}x_0$ and by Proposition \ref{wbI}(1), ${\Downarrow}x_0\subseteq (\mu_X(\alpha))^{\downarrow}$. We will show that
${\Downarrow}x_0\leq (\Phi(r)(\alpha))^{\downarrow}$, so that $r(\mu_X(\alpha)])x_0\leq \bigvee (\Phi(r)(\alpha))^{\downarrow}=r(\Phi(r)(\alpha))$.

If $x\ll x_0$, then there exists $y\in X$ such that $x\ll y\ll x_0$. Since $(\mu_X(\alpha))^{\downarrow}$ is an ideal, by Proposition \ref{wbI}(2), we have $y\in(\mu_X(\alpha))^{\downarrow}$ and thus there exists $A\in\sigma(X)$ such that $\phi(A)\in\alpha$ and $A\subseteq{\uparrow}y$. Then $y\leq\bigwedge A$ and $x\ll\bigwedge A$. By Lemma 2.3, $r^\leftarrow({\Uparrow}x)\in\alpha$. Since ${\Uparrow}x\in\sigma(X)$ and ${\Uparrow}x\subseteq{\uparrow}x$, we have $x\in\Phi(r)(\alpha)$. Hence, ${\Downarrow}x_0\leq (\Phi(r)(\alpha))^{\downarrow}$. This completes the proof of this step.

Therefore, $r\cdot\Phi  r=r\cdot \mu_X$.}

\begin{Pro}$r\cdot\eta_X=id_X$.\end{Pro}

\p{Let $x\in X$. Then by Proposition 1.8, it holds that $r\cdot\eta_X(x)=r([x])=\bigvee[x]^\downarrow=x$.}

\begin{Th} If $(X,r)$ is a $\Phi$-algebra over ${\bf Top}_0$,
then by considering $X$ with the specialization order,
$X$ is a continuous lattice and $r(v)=\bigvee v^{\downarrow}$.\end{Th}

\begin{Lem}For every $A\in\mathcal{O}(X)$, it holds that $A\subseteq r^\rightarrow(\phi(A))$.\end{Lem}

\p{For every $x\in A$, we have $[x]\in\phi(A)$ and then $x=r(\eta_X(x))=r([x])\in r^\rightarrow(\phi(A))$. Hence, $A\subseteq r^\rightarrow(\phi(A))$.}

\begin{Lem}For all $v,w\in\Phi(X)$, if $v\subseteq w$, then $r(v)\leq_{\mathcal{O}(X)}r(w)$.\end{Lem}

\p{Let $A\in\mathcal{O}(X)$ and $r(v)\in A$. Then $v\in r^\leftarrow(A)\in\mathcal{O}(\Phi(X))$. Since $r^\leftarrow(A)$ is an upper set, we have $w\in r^\leftarrow(A)$ and then $r(w)\in A$. By the arbitrariness of $A$, $r(v)\leq_{\mathcal{O}(X)}r(w)$.}

\begin{Pro}The pair $(X,\leq_{\mathcal{O}(X)})$ is a complete lattice, where $\bigwedge A=r([A])$.\end{Pro}

\p{Let $A\subseteq X$. Firstly, if $x\in A$, then $[A]\subseteq [x]$ and then by Lemma 2.8, $r([A])\leq r([x])=x$. That is to say, $r([A])$ is a lower bound of $A$.

Suppose that $y$ is another lower bound of $A$, i.e., $A\subseteq{\uparrow}y$. If $B\in[y]$, then ${\uparrow}y\subseteq B$ and then $B\in[A]$. Thus $[y]\subseteq [A]$ and $y=r([y])\leq r([A])$.

Hence, $\bigwedge A=r([A])$.}

\begin{Lem}Let $\mathcal{A}$ be a directed subset of $\Phi(X)$. Define $\widetilde{\mathcal{A}}\subseteq\mathcal{O}(\Phi(X))$ by $$\widetilde{\mathcal{A}}=\{W\in\mathcal{O}(\Phi(X))|\ \mathcal{A}\cap W\not=\emptyset\}.$$Then $\widetilde{\mathcal{A}}\in\Phi^2(X)$ and
$\mu_X(\widetilde{\mathcal{A}})=\bigvee\mathcal{A}$.\end{Lem}

\p{For all $W_1,W_2\in\mathcal{O}(\Phi (X))$, it is clear that if $W_1\cap W_2\in\widetilde{\mathcal{A}}$ then $W_1,W_2\in\widetilde{\mathcal{A}}$. Conversely, suppose that $W_1,W_2\in\widetilde{\mathcal{A}}$. Then $\mathcal{A}\cap W_1\not=\emptyset$ and $\mathcal{A}\cap W_2\not=\emptyset$, and then there exist $v_1,v_2\in\mathcal{A}$ such that $v_1\in W_1,\ v_2\in W_2$. Since $\mathcal{A}$ is directed, there exists $w\in\mathcal{A}$ such that $v_1,v_2\subseteq w$. We know that $W_1,W_2$ are upper sets of $(\Phi(X),\subseteq)$, we have $w\in W_1\cap W_2$. Thus, $W_1\cap W_2\in\widetilde{\mathcal{A}}$. Hence, $\widetilde{\mathcal{A}}\in\Phi^2(X)$.

For each $A\in\mathcal{O}(X)$, $A\in\mu_X(\widetilde{\mathcal{A}})$ if{}f $\phi(A)\in\widetilde{\mathcal{A}}$ if{}f $\phi(A)\cap\mathcal{A}\not=\emptyset$ if{}f there exists $v\in\mathcal{A}$ such that $A\in u$, if{}f $A\in\bigcup\mathcal{A}$. Hence, $\mu_X(\widetilde{\mathcal{A}})=\bigvee\mathcal{A}$.}

\begin{Pro}$r:\Phi(X)\longrightarrow X$ preserves suprema of directed subsets.\end{Pro}

\p{Let $\mathcal{A}$ be a directed subset of $\Phi(X)$. Then by Lemma 2.10,
$$r(\bigvee\mathcal{A})=r\cdot\mu_X(\widetilde{\mathcal{A}})=r\cdot\Phi(r)(\widetilde{\mathcal{A}}).$$
We will prove that $r\cdot\Phi(r)(\widetilde{\mathcal{A}})\leq\bigvee r_L^\rightarrow(\mathcal{A})$ (the inverse inequality is routine).
Since $\bigvee r^\rightarrow(\mathcal{A})=\bigwedge(r^\rightarrow(\mathcal{A}))^{\uparrow}=r([(r^\rightarrow(\mathcal{A}))^{\uparrow}])$, we only need to show that $\Phi(r)(\widetilde{\mathcal{A}})\subseteq[(r^\rightarrow(\mathcal{A}))^{\uparrow}]$.

In fact, suppose that $A\in\Phi(r)(\widetilde{\mathcal{A}})$. Then $r^\leftarrow(A)\in\widetilde{\mathcal{A}}$ and then $\mathcal{A}\cap r^\leftarrow(A)\not=\emptyset$. Thus there exists $v\in \mathcal{A}$ and $r(v)\in A$.
Suppose that $x\in (r^\rightarrow(\mathcal{A}))^{\uparrow}$. By $v\in \mathcal{A}$, we know that $r(v)\in r^\rightarrow(\mathcal{A})$ and then $r(v)\leq x$, which implies that $x\in A$. Hence, by the arbitrariness of $x$, we have $(r^\rightarrow(\mathcal{A}))^{\uparrow}\subseteq A$; and by the arbitrariness of $A$, we have  $\Phi(r)(\widetilde{\mathcal{A}})\subseteq[(r^\rightarrow(\mathcal{A}))^{\uparrow}]$.
This completes the proof.}

\begin{Lem}For every $v\in\Phi(X)$, define $\mathcal{A}_v\subseteq {\Phi(X)}$ by
$$\mathcal{A}_v=\{[A]|\ A\in v\}.$$ Then $\mathcal{A}_v$ is directed in $\Phi(X)$ and $\bigvee\mathcal{A}_v=v$.\end{Lem}

\p{Firstly, it is obvious that $\mathcal{A}_v\not=\emptyset$. For all $[A], [B]\in\mathcal{A}_v$ with $A,B\in v$, we have $A\cap B\in v$ and then $[A\cap B]\in\mathcal{A}_v$. Clearly, $[A],[B]\subseteq[A\cap B]$
Hence, $\mathcal{A}_v$ is directed. Secondly, $v=\bigcup_{A\in v}[A]=\bigvee\mathcal{A}_v$.}

\begin{Pro}For $v\in\Phi(X)$, $r(v)=\bigvee v^{\downarrow}$.\end{Pro}

\p{We need to show that $r(v)$ is the least upper bound of $v^{\downarrow}$.

Firstly, for each $x\in v^{\downarrow}$, there exists $A\in v$ such that $x\in A^{\downarrow}$. Then $[A]\subseteq v$ and $x\leq\bigwedge A=r([A])\leq r(v)$. That is to say, $r(v)$ is an upper bound of $v^{\downarrow}$.

Secondly, suppose that $y$ is another upper bound of $v^{\downarrow}$. Since $r(v)=r(\bigvee\mathcal{A}_v)=\bigvee r^\rightarrow(\mathcal{A}_v)$, we only need to show that $z\leq y$ for each $z\in r^\rightarrow(\mathcal{A}_v)$. In fact, for $z\in r^\rightarrow(\mathcal{A}_v)$, there exists $A\in v$ such that $z=r([A])=\bigwedge A$. Then $z\in A^\downarrow\subseteq v^\downarrow$ and then $z\leq y$.

Hence, $r(v)$ is the least upper bound of $v^{\downarrow}$, i.e., $r(v)=\bigvee v^{\downarrow}$.}

\begin{Pro}$(X,e_{\mathcal{O}(X)})$ is a continuous lattice.\end{Pro}

\p{By Proposition 2.13, $x=r([x])=\bigvee[x]^{\downarrow}$ holds for every $x\in X$. By Propositions 1.8 and 2.9, $X$ is a continuous lattice.}

\end{document}